\documentstyle{article}

\begin{document}

\newcommand{\uu}{\rza\mapsto\rza}

\def\dj{d\kern-0.4em\char"16\kern-0.1em}

\newcommand{\qed}{\hfill $\Box$ \bigskip}

\newcommand{\dkz}{\noindent {\bf Proof:}{\hspace{0.3cm}} \nopagebreak}
\newcommand{\HDS}{\vrule width0pt height2.3ex depth1.05ex\displaystyle}
\newcommand{\cls}{${\cal C}(\pi_1,\pi_2)$}

\def\teo#1{{\noindent \bf Theorem #1}\hspace{1em}}
\def\lema#1{{\noindent \bf Lemma #1}\hspace{1em}}
\def\de{\noindent {\bf Definition.}\hspace{1em}}
\def\dia{\nopagebreak \hfill$\diamond$}

\def\tekst#1{\mbox{\rm #1}}

\def\lprav{the se\-quen\-ce, $\cal R$, of in\-fe\-rence rules}
\def\dprav{rules from $\cal R$, applied in the same
order as on the left, each of them, possibly se\-ve\-ral times}

\newcommand{\bsim}{\mbox{\boldmath\large$\sim$}}
\newcommand{\bneg}{\mbox{\boldmath\large$\sim$}}
\def\nsim#1{\underbrace{\sim\dots\sim}_{#1\,{\rm times}}}
\def\nneg#1{\underbrace{\sim\dots\sim}_{#1\,{\rm times}}}
\def\nbsim#1{\underbrace{\bsim\dots\bsim}_{#1\,{\rm times}}}
\def\nbneg#1{\underbrace{\bneg\dots\bneg}_{#1\,{\rm times}}}

\def\o#1{\overline{#1}}

\def\dz{-^\ast}
\def\lz{^\ast-}

\def\lz{\langle}
\def\rz{\rangle}

\def\dzi{{\rm i}^\ast}
\def\dze{{\rm e}^\ast}

\def\k#1#2{\foramebox[#1 em][l]{$#2$}}

\def\aks#1{#1}
\def\f#1#2{{{\HDS #1}\over{\HDS #2}}}
\def\fp#1#2{{{\HDS #1}\atop{\HDS #2}}}
\def\afrac#1{{\phantom{\HDS #1}\atop{\HDS #1}}}
\def\bfrac#1#2{{\atop{\phantom{
\vrule width0pt height2.3ex depth#2\displaystyle \Theta}\atop{\HDS
#1}}}}

\def\rza{{\mbox{\hspace{0.5em}}}}
\def\rzb{{\mbox{\hspace{2em}}}}
\def\rzbc{{\mbox{\hspace{4,5em}}}}
\def\rzc{{\mbox{\hspace{7em}}}}
\def\rzd{{\mbox{\hspace{11em}}}}

\def\Bpak#1#2#3{\f{\f{#1}{\cdots} \pravila{#2}}{#3}}
\def\Epak#1#2#3{\f{\f{#1}{\cdots} \Epravila{#2}}{#3}}

\def\BBpak#1#2#3{\f{#1}{\f{\cdots}{#3}} \Bpravila{#2}}
\def\Cpak#1#2#3{\f{\f{#1}{\cdots} \ppravila{#2}}{#3}}
\def\Dpak#1#2#3{\f{\f{#1}{\cdots} \pppravila{#2}}{#3}}

\def\p#1{\{#1\}}

\def\prj#1{\enskip{\scriptstyle(\rm{#1})}}
\def\prjm#1{\enskip{\scriptstyle({#1})}}
\def\pr#1#2{\enskip{\scriptstyle({#1}\enskip\rm{#2})}}
\def\prav#1#2{ \rza \makebox[-.5em][l]{\mbox{\scriptsize{($#1$\enskip\rm #2)}}}}
\def\pravj#1{ \rza \makebox[-.5em][l]{\mbox{\scriptsize{(\rm #1)}}}}
\def\pravjm#1{ \rza \makebox[-.5em][l]{\scriptsize{($#1$)}}}

\def\fns#1{{\footnotesize{#1}}}
\def\ss#1{{\scriptsize{#1}}}

\def\pravilo#1{ \rza \makebox[-.5em][l]{\mbox{\scriptsize{\rm #1}}}}
\def\lpravilo#1{ \makebox[-.5em][r]{\mbox{\rm #1}} \rza}
\def\pravila#1{ \rza \makebox[-.5em][l]{\raisebox{-1.7ex}
{\parbox{15ex}{\footnotesize \baselineskip=0.5\baselineskip #1}}}}

\def\Epravila#1{ \rza \makebox[-.5em][l]{\raisebox{-1.7ex}
{\parbox{25ex}{\footnotesize \baselineskip=0.5\baselineskip #1}}}}

\def\Bpravila#1{ \rza \makebox[-.5em][l]{\raisebox{-1.7ex}
{\parbox{30ex}{\footnotesize \baselineskip=0.5\baselineskip #1}}}}

\def\ppravila#1{ \rza \makebox[-.5em][l]{\raisebox{-1.7ex}
{\parbox{19ex}{\footnotesize \baselineskip=0.5\baselineskip #1}}}}
\def\pppravila#1{ \rza \makebox[-.5em][l]{\raisebox{-1.7ex}
{\parbox{30ex}{\footnotesize \baselineskip=0.5\baselineskip #1}}}}
\def\naz#1{$\tekst{#1}^z$}
\def\pravilo#1{ \rza \makebox[-.5em][l]{\mbox{\rm #1}}}
\def\lpravilo#1{ \makebox[-.5em][r]{\mbox{\rm #1}} \rza}
\def\pravila#1{ \rza \makebox[-.5em][l]{\raisebox{-1.7ex}
{\parbox{15ex}{\footnotesize \baselineskip=0.5\baselineskip #1}}}}
\def\ppravila#1{ \rza \makebox[-.5em][l]{\raisebox{-1.7ex}
{\parbox{19ex}{\footnotesize \baselineskip=0.5\baselineskip #1}}}}
\def\pppravila#1{ \rza \makebox[-.5em][l]{\raisebox{-1.7ex}
{\parbox{30ex}{\footnotesize \baselineskip=0.5\baselineskip #1}}}}
\def\naz#1{$\tekst{#1}^z$}

\def \ra{\rightarrow}
\def\lra{\leftrightarrow}
\def\i{\wedge}
\def\ili{\vee}
\def\xa{\alpha}
\def\xb{\beta}
\def\xg{\gamma}
\def\xd{\delta}
\def\Xg{\Gamma}
\def\Xd{\Delta}
\def\Xs{\Sigma}
\def\Xl{\Lambda}
\def\Xt{\Theta}
\def\XI{\rm I}
\def\XE{\rm E}
\def\XV{\rm V}
\def\ol{\overline}
\def\ul{\underline}
\def\ua{\uparrow}
\def\da{\downarrow}
\def\sr{\rm sr}
\def\com{\rm com}
\def\max{\rm max}
\def\er{\rm e}
\def\fr{\rm fr}
\def\fl{\rm fl}
\def\Int{\rm Int}


\noindent DISCUSSION PAPER

\vskip6truemm \centerline{
\bf Arrow--Sen theory simplified\rm}
\vskip2truemm

\centerline{\it Branislav Bori\v ci\'c}
\centerline{\it University of Belgrade, Faculty of Economics and Business}
\centerline{\it boricic@ekof.bg.ac.rs}
\vskip3truemm

\bf Abstract. \rm The traditional Arrow--Sen Social Choice Theory $\bf{TSCT}$ is a mathematical theory built apparently on higher--order formal language. In this paper, we propose a reformulation and reclassification of the $\bf{TSCT}$ axioms in order to obtain a simpler theory based on the first--order language axioms, keeping the spirit of original ideas. This new theory, called Simplified Social Choice Theory, denoted by $\bf{SSCT}$, presents a sub--theory of $\bf{TSCT}$. Roughly speaking, we extract all quatifications over $n$--tuples of binary relations from the axioms of $\bf{TSCT}$ and move them to the meta--level obtaining a sub--theory $\bf{SSCT}$ of $\bf{TSCT}$. More accurately, we assign to each traditional higher--order axiom $\bf{TA}$ its simplified first--order version $\bf{SA}$ such that $\bf{TA}\vdash\bf{SA}$, i.e. $\bf{SA}$ can be logically derived from $\bf{TA}$. In this way we define a simpler and more accessible set of principles that provide a context in which we can prove many propositions analogous to well--known theorems including the Arrow's impossibility of Paretian non--dictatorship and Sen's impossibility of Paretian liberal. These simplifications are the result of decades of lecturing by the author with the aim of bridging the barriers between this beautiful complex theory and his students.\vskip3truemm



\vskip3truemm
\S 1. \bf Introduction. \rm The teaching process requires the introduction of the new concepts, step by step, in order to be understandable and clear to most of our students. In presenting the elements of Social Choice Theory, the first usual step is to sketch out the basics of preference logic, i.e. to define and motivate strict, weak and indifference preference relations. The axioms concerning preferences are simple and contain quantification over alternatives only. For instance, the transitivity axiom for a binary relation $R$ has the following form:
$$(\forall x,y,z\in X)(xRy\wedge yRz\to xRz)$$
where $X$ is a fixed finite set of all possible alternatives, and this formula, well known to students from other parts of mathematics, evidently, belongs to the first--order language. The next step is devoted to explaining the motives for the differences and congruences between the relationship of individual and social preferences, which is not so much a mathematical but a philosophical and sociopsychological controversy.

Concepts such as a social welfare function, including quantification over relations or over profiles, seem to be the first stage that is too difficult for students.

During my personal efforts to find, for my classes, more understandable proofs and forms of Arrow's Impossibility Theorem, I found that there are at least four barriers to be overcome. The first one was the methodological controversy of impossibility (see [BBa]). The second one was the logical formalism of preference logic. The third were psychological and philosophical ambiguities regarding two kinds of preferences, individual and social ones. Finally, the fourth and the strongest barrier was the complex mathematical context in which the Theorem was proved. All those where motives for looking for easier formulations, proofs, context, and examples.

After forty years of the author's teaching experience with the students of economics, business and statistics, who only have a limited  interest in abstract theoretical mathematical themes, but with great enthusiasm for practical and applicable knowledge, my teaching career could be considered a fight with difficult complexities and dangerous simplifications, with defining practical procedures justified by theoretical argumentation, for a sole purpose of fine tuning of transfer, in an adequate way, of knowledge acceptable and useful for students. In this broad context, this paper presents an attempt to separate a set of logically and mathematically simpler axioms of traditional Arrow--Sen Social Choice Theory and to enable better understanding the spirit of crucial results of this complicated theory through an easier approach.

First of all, we tend to replace the usual pure mathematical argumentation and presentation of impossibility theorems in Social Choice Theory by a logical and methodological one, based on terms of consistency and inconsistency, which seems to be more intuitive, natural and understandable. In particular, we propose a method of simplification of some axioms of traditional Arrow--Sen's Social Choice Theory based on elementary logical argumentation. Then, we try to discuss the relationships between practice and theory, in general, but also in the particular field of social decision making.

As a common point for methodology, practice, applied mathematics, logic and education, Arrow's theorem appears in an informal descriptive and simplified form even in basic textbooks on economics, but we are not sure that its essential idea presented in this way is understandable enough. Any proof of inconsistency needs a very formal context. This is the reason why Arrow--Sen theory is based on a high level formalism thus presenting by itself a barrier to an easy understanding of the whole concept. The importance of Arrow's and Sen's results encouraged numerous attempts to explain them, including this paper and author's efforts to make impossibility results approachable to students.

In this paper, in order to obtain a simpler description of Social Choice Theory, based on inevitable formalism, we try to substitute the notion of impossibility by deducibility (or inconsistency), as a more logical concept, to avoid the notions of social welfare function and social decision function, and to simplify the form of a social choice axioms in which the quantifications over alternatives, over individuals and over profiles appear simultaneously. In that process, we get some simpler fragments of traditional Social Choice Theory allowing us to express some situations that resemble some of the known impossibility results and to raise questions about their provability in this new context.

In the last century Arrow's and Sen's works\footnote{\rm K. J. Arrow (1921--2017), Nobel prize in economics 1972; A. K. Sen (1933--), Nobel prize in economics 1998.}, based on an axiomatic method, promoted the most influential ideas in Social Choice Theory. This theory (see [AR], [AS] or [ASS]) opened a list of vaious problems including the ethical and philosophical dimension of basic requirements for a group decision processes (see [BUCa], [BUCb], [Sc], [Se] or [SZ]), as well as, discussion regarding the degree of acceptable formalism in defining these processes (see [Sa], [Se], [RTH], [MP], [NIP], [THW], [AHW], 
[PER], [PY]
or [POR]). One of particularly stimulating circumstance was that Arrow's impossibility theorem testified an impossibility and the numerous papers published later were aimed to redefine the formal normative framework ensuring the possibility for making decisions (see [Sa]).

One of the first significant logical analysis of Arrow's impossibility theorem was made by R. Routhley in [RTH]. Over the last decade, there have been different approaches to analyzing the complexity of various aspects of social choice theory. The first step is always to find a simpler language in which the basic statements of the theory can be expressed and interpreted. A common bond connecting those papers is that all of them try to illuminate difficult facts and argumentations of social choice theory from different aspects. M. Pauly [MP] proposed to compare formalizations in social choice theory according to the language used for expressing the axioms and suggests a formalist approach to axiomatization results which uses a restricted formal logical language to express axioms. The advantages of his approach include the possibility of non--axiomatizability results, a distinction between absolute and relative axiomatizations, and the possibility to ask how rich a language needs to be to express certain axioms, arguing for favoring axiomatizations in the weakest language possible. We see that our approach is consistent with the basic ideas presented in [MP]. T. Nipkow [NIP] discussed proofs of Arrow's impossibility theorem formalized in a higher--order language. In [TL], P. Tang and F. Lin proposed alternative proof of Arrow's and other impossibility theorems with the idea to use induction in order to reduce the theorems to the base case with $3$ alternatives and $2$ agents and then use computers to verify the base case. N. Troquard, W. van der Hoek and M. J. Wooldridge in [THW] introduced a logic specifically designed to support reasoning about social choice functions enabling to show that every social choice function can be characterised as a formula of this logic. Also, T. Agotnes, W. van der Hoek and M. J. Wooldridge in [AHW] presented a sound and complete axiomatization of a modal logic supporting reasoning about judgment aggregation scenarios and preference aggregation, where the logical language is interpreted directly in judgment aggregation rules. 
T. Perkov [PER] introduced a sound and complete natural deduction system for modal logic of judgment aggregation making it possible to express classical properties of judgment aggregation rules and famous results of social choice theory, including Arrow's impossibility theorem. In [PY], E. Pacuit and F. Yang developed a version of independence logic that can express Arrow's properties of preference aggregation functions, and then they proved that Arrow's impossibility theorem is derivable in a natural deduction system for the first--order consequences of their logic. 
D. Porello (see [POR]) introduced a number of logics, based on a substructural propositional logic that allows for circumventing inconsistent outcomes, for modelling collective propositional attitudes that are defined by means of the majority rule.

Above mentioned papers, dealing with formalization of social choice theory, can be roughly divided into two sorts: those working in (first-- or higher--order) predicate languages, e.g. [RTH], [MP], [NIP], [TL], 
[THW] and [PY], 
and those based on modal languages, e.g. [AHW], [PER] and [POR]. Unlike them, we consider and justify some simplified subtheories of social choice theory, i.e. its pure first--order fragments enabling us to explain some impossibilities. The next step we take in this simplification process is treating these impossibilities in a propositional language extended by a finite list of binary (preference) relations defined over propositional formulae as alternatives. In this step we rely on concept of G. H. von Wright's preference logic (see [VWa] and [VWb]) and idea of combining logics (see [GB]) as given in our earlier works [BBd] and [BBe].

We divide the traditional social choice axioms introduced by K. Arrow and A. Sen into two classes, based on their linguistic and mathematical complexity. The first class consists of 'the unrestricted domain' $\bf{U}$, and 'the independence of irrelevant alternatives' $\bf{IIA}$, which, both of them, including some modified versions of $\bf{IIA}$ (see [MSK]), need a language of higher--order (see also [FISH] and [HP]), and which can be treated as a kind of meta--axioms. The second class contains a group of linguistically simpler axioms, such as dictatorship, weak dictatorship, vetoing, liberalism and the Pareto rule. Naturally, it is possible to make an easier logical analysis of deductive properties and relationships between axioms belonging to the second class.

We think it is natural that individual preferences have an ordinal (but not cardinal) character, due to the fact that individuals usually have vague opinions about most alternatives. For instance, in political life individuals are commonly able to rank just a few parties, and most of them left unranked, so that the individual preferences are even incomplete (non--linear) i.e. there are alternatives which some individuals cannot rank. On the other side, in practice, social preference resulting from individual ones are quite natural to be complete (linear) and even numerically ranked. For instance, parliamentary elections generate a linear social numerical ranking from preferences presenting non--linear non--numerical descriptive individual rankings. Namely, voters participating in political elections usually do not have any attitude about all parties or political groups taking parts in the campaign, but only about a few of them and those attitudes are just on a rough descriptive level, as 'good', 'better than' or 'bad'. In contrast to the vague and incomplete individual preferences, the results of elections are commonly expressed as complete numerical ranking of all participants in elections through the percentage of votes or final number of seats in parliament. A possible conclusion which can be drawn is that, in practice, individuals do not satisfy Arrow's rational choice axioms (e.g. linearity), but groups and societies do satisfy.

Foundational statements of possible--impossible type define the demarcation line between a domain in which decision procedures give some results and an area where decision procedure does not exist at all. A. Sen shows [Sa] how this line is hardly visible and sensitive, and how it depends on small and seemingly insignificant changes of conditions. Sen introduced the weakening and strengthening of definitions of decision making procedures (a social welfare function, a social decision function) and showed how to obtain some possibility results and 'improvements' from original Arrow's impossibility theorem. In this way the demarcation line between 'a possible area' and 'an impossible one' becomes more visible and it gives explicit instructions on how, in practice, some basic properties in defining a decision making process must not or could be combined. This is of great importance because inconsistencies caused by some conditions are not trivial, and avoiding them helps to define better indisputable voting laws or group decision making rules.

One of those paradoxical conditions which surprisingly causes logical inconsistencies is the Pareto rule. In this paper we complete a list of the Pareto rule's paradoxes starting with Arrow's impossibility theorem, which can be understood as 'the impossibility of a Paretian non--dictator', and Sen's 'the impossibility of a Paretian liberal', concluding with Mas--Colell's and Sonnenschein's 'the impossibility of a Paretian non--vetoer' and our 'impossibility of a non--Paretian dictator'. All these facts can be considered as a list of 'Paretian paradoxes'.

J. M. Buchanan's\footnote{\rm J. M. Buchanan (1919--2013), Nobel prize in economics 1986.} criticism seems to be still actual (see [BUCa], [BUCb] or [Se]). The Pareto rule, presenting the most controversial condition from the ethical stand point, becomes the most problematic axiom from the logical point of view as well. Ethically 'better to all', comprehended as 'better to society', opens the possibility for a small improvement for a large circle of society members and a big improvement for a small circle of society members, implying remarkable growth of social inequalities, which cannot characterized as a better society state. On the other side, logically, the presence of the Pareto rule, as we see, provokes inconsistencies.

M. Balinski and R. Laraki [BL] leave the concept of pure ordinal preferences by insisting on evaluation of the alternatives instead of defining preferences over the alternatives. Their concept, based on the so called majority judgement is, of course, more applicable and effective than the pure ordinal one, but there is a question if individuals (voters) are competent to evaluate the alternatives independently when evaluations have a descriptive ('excellent', 'very good',...) or, numerical ($5$, $4$,...) character. H. de Swart [DS] includes in his brilliant book a modern approachable introductory chapter devoted to voting systems, with elements of ordinal Social Choice Theory, where the idea of majority judgement voting rule is discussed in details.

This paper is organized as follows.

First we give a brief review of Arrow--Sen social choice theory, relevant to our approach, where we try to define a demarcation between 'complex' and 'simple' axioms, founded on the expressibility of axioms in the first--order predicate language. Axioms of 'unrestricted domain' $\bf{U}$ and 'independence of irrelevant alternatives' $\bf{IIA}$, belonging to the class of 'complex' conditions, were introduced by Arrow (see [AR], [FISH], [KEL], [Sc], [SCH], [TAY]). Our list of 'simple' axioms starts with Arrow's traditional 'dictatorship' $\bf{TD}$ and 'the Pareto property' $\bf{TP}$. Sen introduced the axiom of 'liberalism' $\bf{TL}$, during his further development of the theory, and pointed out to the 'liberal paradox' enriching the theory with a new approachable example of impossibility (see [Sb], [Sc], [Sd], [Se], [BBb], [BBd], [BBe], [BBf]). The vetoer axiom $\bf{TV}$, introduced by Mas--Colell and Sonnenschein [MS], was considered by Fishburn [FISH] as part of an ordinal treatment, characteristic for the Arrow--Sen theory, but the status of the veto power in weighted voting systems was discussed in [TAY] and [WIN], and in many other works (see [BS] and [SRE]).

In the sequel, we propose and justify a method of simplification of traditional axioms. This method reduces parts of Social Choice Theory to some of its fragments based on the first--order language, but rich enough to deal with their interdeducibilities and counterparts of traditional impossibilities. Roughly speaking, the method consists of accepting the forms $\exists xA$ instead of $\exists x\forall yA$, based on general logical law $\exists x\forall yA\to\forall y\exists xA$ and then moving universal quantification on meta--level.

This method of simplification produces a new list of simplified axioms dictatorship $\bf{SD}$, the Pareto rule $\bf{SP}$, liberalism $\bf{SL}$ and vetoer $\bf{SV}$, expressed in first--order language. In this new context we establish some expected deductive interdependence and prove the counterparts of well--known statements of social choice theory, including Arrow's and Sen's impossibility theorems.

\bf Highlights. \rm Inspired by Traditional Social Choice Theory $\bf{TSCT}$, prefixed by $\bf{T}$ --- for 'traditional', dealing with axioms such as dictatorship $\bf{TD}$, liberalism $\bf{TL}$, vetoing $\bf{TV}$ and the Pareto rule $\bf{TP}$, but also supposing the presence of 'unrestricted domain' $\bf{U}$ and 'the independence of irrelevant alternatives' $\bf{IIA}$, we consider a fragment of so--called Simplified Social Choice Theory $\bf{SSCT}$, prefixed by $\bf{S}$ --- for 'simplified', based on new simplified axioms of dictatorship $\bf{SD}$, liberalism $\bf{SL}$, vetoing $\bf{SV}$ and the Pareto rule $\bf{SP}$, expressed exclusively in first--order language of predicates, such that each traditional axiom \it deductively implies \rm its simplified version:
$$\bf{TD}\vdash\bf{SD}, \ \ \bf{TL}\vdash\bf{SL}, \ \ \bf{TV}\vdash\bf{SV}\ \ {\rm and} \ \ \bf{TP}\vdash\bf{SP}$$
where, we also suppose that conditions $\bf{U}$ and $\bf{IIA}$ hold in $\bf{SSCT}$ in a specific way on a metatheoretical level.

Our aim is to define a new context, $\bf{SSCT}$, to analyze pure logical relationships between these four axioms in order to obtain counterparts to some known impossibility theorems, by avoiding complex mathematical machinery such as social welfare function and so on. Bearing in mind deductive connections between $\bf{T}$--axioms and $\bf{S}$--axioms, the $\bf{SSCT}$ presents a subtheory of $\bf{TSCT}$, but this fact does not indicate that each impossibility provable in $\bf{SSCT}$ \it a fortiori \rm implies the corresponding impossibility in $\bf{TSCT}$. For instance, neither counterpart $\bf{SP}\vdash\bf{SD}$ implies original Arrow's theorem $\bf{TP}\vdash\bf{TD}$, nor vice versa. Consequently, these two statements can be considered as two roughly connected facts in two parallel worlds. On the other side, from $\bf{TP}\vdash\bf{SP}$ and $\bf{TL}\vdash\bf{SL}$, we can directly derive well--known Sen's 'impossibility of a Paretian liberal': $\bf{TP},\bf{TL}\vdash$, from its simplified version $\bf{SP},\bf{SL}\vdash$, meaning that the axioms $\bf{SP}$ and $\bf{SL}$, and, consequently, the axioms $\bf{TP}$ and $\bf{TL}$, when appear together, make a theory inconsistent.

We consider that the value of this simplified approach is in giving an opportunity to a wider circle of readers to understand the basic ideas, results and spirit of traditional Social Choice Theory better.

Finally, let us describe our framework more precisely: (i) 'Unrestricted domain' $\bf{U}$, and 'the independence of irrelevant alternatives' $\bf{IIA}$, are assumed to be general meta--conditions of our theory; (ii) instead of traditional dictatorship $\bf{TD}$, liberalism $\bf{TL}$, vetoing $\bf{TV}$ and the Pareto rule $\bf{TP}$, we use their weaker simplified versions: $\bf{SD}$, $\bf{SL}$, $\bf{SV}$ and $\bf{SP}$, respectively; (iii) in our theory we are able to formulate and prove counterparts of known impossibility results of $\bf{TSCT}$.

\vskip3truemm
\S 2. \bf Traditional Arrow--Sen Theory --- An Outline of Traditional Approach. \rm Arrow--Sen theory covers both macro level of social choice, such as a parliamentary or presidential election procedure in a country, and micro level of group decision making, like consumer preferences, statutory voting or deciding in a joint stock company.

The original proof of Arrow's Impossibility Theorem was approachable and understandable to a very limited circle of scientists during the first two decades after its publishing.

Arrow--Sen approach to Social Choice Theory can be considered a mathematical (and logical) analysis of conditions under which a social decision function exists. Any procedure generating a general social opinion from particular individual opinions, presents a social decision function. Our approach is essentially based on traditional books concerning this subject (see [AR], [KEL], [Sc], [SCH] or [TAY]).

More accurately, this procedure deals with social choices over a finite set $X$ of alternatives, respecting preferences of individuals of a finite set $V$. It is supposed that individuals and society satisfy rational choice axioms, i.e. that each individual preference relation $R_i$, characterizing the behavior of an individual $i\in V$, is linear, $(\forall x,y\in X)(xR_iy\vee yR_ix)$, and transitive, $(\forall x,y,z\in X)(xR_iy\wedge yR_iz\to xR_iz)$, and that the corresponding preference relation $R$, characterizing the behavior of society, is also linear and transitive. Such relations $R_i$ and $R$ are called the weak preference relations, and $xRy$ stands for '$x$ being regarded as at least as good as $y$'. Each weak preference relation $R$ defines corresponding strict preference $P$, and indifference $I$, as follows: $xPy$ iff $xRy\wedge\neg yRx$, and $xIy$ iff $xRy\wedge yRx$. In this case we also have that an indifference relation $I$ can be generated by a strict preference relation $P$: $xIy$ iff $\neg(xPy)\wedge\neg(yPx)$, as well as its weak version $R$: $xRy$ iff $xPy\vee xIy$.

Here we use symbols for universal, $\forall$, and existential, $\exists$, quantifiers, as well as, propositional connectives for negation, $\neg$, conjunction, $\wedge$, disjunction, $\vee$, implication, $\to$, and equivalence, $\leftrightarrow$, with the usual meaning they have in classical logic. Also, we use the turnstile symbol, $\vdash$, for deduction relation in an informal way, $A,B\vdash C$, in order to express that "$C$ can be derived from $A$ and $B$".

As we can see, a wider context for Arrow--Sen theory is the classical set theory, more descriptive than the formal one.

This relation $R$ enables us to define a choice set $C(Y,R)=\{ x|x\in Y\wedge(\forall y\in Y)xRy\}$, presenting the set of best alternatives, with respect to $R$ and $Y\subseteq X$, while $C(Y,R)$ will present a choice function, over the set of all alternatives, if $C(Y,R)$ is non--empty for every non--empty $Y\subseteq X$. According to [Sa], a rule is defined as a functional relation specifying one and only one social binary relation $R$ for each profile of individual ordering $(R_1,\dots,R_n)$, with one $R_i$ for each individual $i\in V$ $(1\leq i\leq n)$. A social welfare function is defined as a rule the range of which is restricted to the set of orderings. A social decision function is defined as a rule ranged over relations $R$ generating a choice function $C(Y,R)$ over entire $X$.

Arrow's original theory is founded on the following four axioms:

Axiom $\bf{U}$ of 'unrestricted domain' requires that the choice function can be applied to {\bf any profile} of logically possible profiles of individual preferences.

Axiom $\bf{IIA}$ of 'the independence of irrelevant alternatives' ensures that for any binary relations $R$ and $R^\prime$ generated respectively by {\bf any two profiles}, $n$--tuples of individual preferences $(R_1,\dots,R_n)$ and $(R_1^\prime,\dots,R_n^\prime)$, and for all pairs of alternatives $(x,y)\in Y^2$, where $Y$ is any subset of $X$, if $(\forall i\in V)(xR_iy\leftrightarrow xR_i^{\prime}y)$, then $C(Y,R)=C(Y,R^{\prime})$.

Non--dictatorship axiom $\bf{TND}$ states that there is no person $i\in V$, a dictator, having such power that, for {\bf all profiles} and each two alternatives $x$ and $y$, if $i$ prefers $x$ to $y$, society must prefer $x$ to $y$ as well.

The Pareto property $\bf{TP}$ claims that, for {\bf all profiles}, if every individual $i\in V$ prefers $x$ to $y$, then society must prefer $x$ to $y$. This is, in fact, a weak version of the Pareto principle, as introduced by Arrow (see [AR], [Sa] and [Sb]).

In denotation of some axioms we use the prefix $\bf{T}$ pointing out that this is a 'traditional' form of axiom, e.g. $\bf{TND}$ and $\bf{TP}$. It is necessary because we will exploit the differences between a 'traditional' and 'simplified' form of axioms.

Now, Arrow's famous impossibility theorem, or as it is originally called 'General Possibility Theorem', can be presented as follows (see [AR] and [Sd]):

\vskip3truemm
\bf Arrow's Impossibility Theorem. \it There is no social welfare function satisfying axioms ${\bf{U}}$, ${\bf{IIA}}$, ${\bf{TND}}$ and ${\bf{TP}}$. \rm\vskip3truemm

Sen has taken into consideration, following the spirit of J. S. Mill's liberalism comprehension, 'the liberalism axiom' $\bf{TL}$ (see [Sb], [Sc], [Sd], [BBb], [BBd] or [BBe]) by which, for {\bf all profiles} and each individual $i\in V$ there is at least one pair of alternatives $(x,y)\in X^2$ such that $x\neq y\wedge(xP_iy\to xPy)\wedge(yP_ix\to yPx)$.

Sen's famous result (see [Sb], [Sc] or [Sd]) known as the 'impossibility of a Paretian liberal' or the 'liberal paradox' can be formulated in the following way:

\vskip3truemm
\bf Sen's Impossibility Theorem. \it There is no social decision function satisfying axioms ${\bf{U}}$, ${\bf{TL}}$ and ${\bf{TP}}$. \rm\vskip3truemm

Let us emphasize that Sen in [Sd] makes a subtle difference between a social welfare function and a social decision function. Namely, a social welfare function is ranged over the set of orderings, while a social decision function is ranged only over those binary relations $R$ each of which generates the choice function $C(X,R)$ over entire $Y$. But such details will not be essential for our approach given in the sequel of this paper.

A logical analysis of Arrow--Sen axioms (see [RTH]) shows that they can be divided naturally, on the basis of their complexity, into two groups: $\{{\bf{U}},{\bf{IIA}}\}$ and $\{{\bf{TND}},{\bf{TP}},{\bf{TL}}\}$. The first two axioms, ${\bf{U}}$ and ${\bf{IIA}}$, have a deeply schematic metatheoretical character and we will use them as general properties of our system. On the other side, the statements, such as ${\bf{TND}}$, ${\bf{TP}}$, ${\bf{TL}}$, and their negations, present formally simpler and mutually similar structures enabling an easier logical analysis of their deductive interdependencies.

\vskip3truemm
\S 3. \bf Logical Tools. \rm In order to consider and formulate logical relationships between our statements more formally, we use the turnstile symbol $\vdash$, denoting the deduction relation. Namely, we write $\bf{X}\vdash\bf{Y}$ to express the fact that statement $\bf{Y}$ can be inferred logically from the statement $\bf{X}$. In particular, $\bf{X},\bf{Y}\vdash$ is used to denote that a set $\{\bf{X},\bf{Y}\}$, consisting of statements $\bf{X}$ and $\bf{Y}$, is inconsistent, i.e. that the simultaneous satisfaction of both statements $\bf{X}$ and $\bf{Y}$ is impossible.

Here we want to point out some elementary logical properties of deduction relation. For instance, the derivation of $\bf{X},\bf{Y}\vdash\bf{Z}$, from hypotheses $\bf{Y}\vdash\bf{W}$ and $\bf{X},\bf{W}\vdash\bf{Z}$, is known as \it the cut rule, \rm or, alternatively, as \it the hypothetical syllogism rule. \rm Also, the equiderivability of $\bf{X},\bf{Y}\vdash\bf{Z}$, with both $\bf{X},\bf{Y},\bf{NZ}\vdash$ and $\bf{X},\bf{NZ}\vdash\bf{NY}$, where $\bf{NX}$ denotes the negation of statement $\bf{X}$, follows immediately from the basic logical properties of negation connective and deduction relation (see [DS], [TAK], [DvD] or [BBc]). Let us note that we use $\bf{NX}$, $\bf{X}$ prefixed by $\bf{N}$, instead of traditional logical denotation $\neg\bf{X}$, bearing in mind the original symbolism usually used in social choice theory, where, for instance, ${\bf{D}}$ and ${\bf{ND}}$ denote dictatorship and non--dictatorship conditions, respectively.

This approach makes it possible to formulate impossibility results more easily and formally manipulate with them: from $\bf{X}\vdash\bf{Y}$ we conclude $\bf{X},\bf{NY}\vdash$, i.e. that the theory containing $\bf{X}$ and $\bf{NY}$ as its statements (axioms or inferred conclusions) is not consistent, or, in other words, that these statements make this theory impossible. Also, if a set $\{\bf{X},\bf{Y}\}$ is inconsistent, i.e. $\bf{X},\bf{Y}\vdash$, then $\bf{X}\vdash\bf{NY}$, and vice versa: from $\bf{X}\vdash\bf{NY}$ we can conclude that the theory containing $\bf{X}$ and $\bf{Y}$ is inconsistent.

\vskip3truemm
\S 4. \bf Simplification Method. \rm From a logical view point, the Social Choice Theory can be considered an analysis of deductive interdependence between various (groups of) axioms appearing in a social choice context. This deductive analysis covers at least the two following phenomena: impossibility and complexity. An interesting aspect of each theory is connected with its (in)consistency ((im)possibility), and we know for many such examples in social choice (see [BBa], [BBb], [FISH], [KEL], [Sa], [Sb], [Sc]). On the other hand, the complexity of a theory defines the limits in its understanding and applicability. Our experience says that teaching elements of Social Choice Theory to students who are not mathematicians is a great challenge because of their complexity. Namely, the complexity of Social Choice Theory presents a difficult barrier for reasonable teaching this subject. We try to simplify some things, but every procedure of simplification brings the danger of banalization.

While Routhley's 'repairing proofs' (see [RTH]) focuses on the subtle differences between two logical forms
$$\forall xA(x)\vdash\forall xB(x){\rm \ \ and \ \ \ }A(x)\vdash B(x)$$
and their roles in presenting a proof of Arrow's impossibility theorem, we deal with the problem of defining some fragments of Social Choice Theory which are not too complex and which can be obtained by substituting some traditional axioms of the form
$$\exists x\forall yA$$
by simpler ones
$$\exists xA$$
relying on the general logical fact that
$$\exists x\forall yA\vdash\forall y\exists xA$$
and then moving the universal quantification $\forall y$ to some kind of metatheoretical level. This operation can be of great importance when the object '$\forall y$' belongs essentially to the higher--order language. By this procedure we can obtain a similar but essentially simpler fragment of the theory which could be more approachable than the original one.

In case of simple axioms of traditional Social Choice Theory, such as, for instance, the Pareto rule, dictatorship, vetoer and liberalism, there are two typical quantifier prefixes in the traditional approach:

(1) there exists an individual $i$, such that for {\bf all profiles} $\mathcal{P}$ and all alternatives $x$ and $y$, $A$; and

(2) for {\bf all profiles} $\mathcal{P}$ and all alternatives $x$ and $y$, $A$,

\noindent which, respectively, can be expressed symbolically as
$$\exists i\forall{\mathcal{P}}\forall x\forall yA{\rm \ \ \ and \ \ \ }\forall{\mathcal{P}}\forall x\forall yA$$
and we propose to substitute the first form by its consequence
$$\forall{\mathcal{P}}\exists i\forall x\forall yA$$
in order to get the universal quantification over profiles as a prefix, which does not belong to the first--order language. In this way, we obtain a uniform form for all these axioms, each of them prefixed by the universal quantification over {\bf all profiles}. In the next iteration, we move this quantification over profiles into the level of metatheory, and then we work with fragments or, say, subtheories based on simple first--order axioms.

In short, roughly, from '$\exists i\forall{\mathcal{P}}\forall x\forall yA$' we proceed to '$\forall{\mathcal{P}}\exists i\forall x\forall yA$', and finally 'we suppose that, for {\bf all profiles} $\mathcal{P}$, we have an axiom: $\exists i\forall x\forall yA$'. Let us point out that we consider a schematic character of 'for {\bf all profiles}' similarly as the status of axiom--schemata in propositional logic, where, for instance, instead of a formula $\forall A\forall B(A\to(B\to A))$ we use the following one $A\to(B\to A)$ treated as an axiom--scheme, with 'for {\bf all formulae} $A$ and $B$' moved on meta level, enabling us to avoid a substitution rule as a primitive rule of a system.

This is a way to define an essential simplification of some parts of traditional Social Choice Theory, but we still hope that this simplification has preserved the basic spirit of traditional Social Choice Theory. Moreover, we believe that this approach is logically quite justified and that the corpus of impossibility results can be represented correctly in this way, due to the obvious fact that if a subtheory is inconsistent, then, {\it a fortiori}, each its extension will be inconsistent as well.

It might seem that in this way we have desecrated the authentic Arrow--Sen tradition, but we believe that: (1) we preserved the spirit of formal logical treatment of descriptive social choice conditions; (2) we reduced the complexity degree of conditions under consideration, and (3) we enabled simplified argumentations for some social choice impossibilities.

More accurately, instead of traditional dictatorship
$${\bf{TD}}: \ \ \ \ (\exists i\in V)(\forall{\mathcal{P}})(\forall x,y\in X)(xP_iy\to xPy)$$
and the traditional Pareto rule
$${\bf{TP}}: \ \ \ \ (\forall{\mathcal{P}})(\forall x,y\in X)((\forall i\in V)xP_iy\to xPy)$$
we accept their variations
$${\bf{SD}}: \ \ \ \ (\exists i\in V)(\forall x,y\in X)(xP_iy\to xPy)$$
and
$${\bf{SP}}: \ \ \ \ (\forall x,y\in X)((\forall i\in V)xP_iy\to xPy)$$
supposing that these variations hold for {\bf all profiles} $\mathcal{P}$, which is in line with the general assumption about the schematic character of axioms. We emphasize that in both cases we have:
$$\bf{TD}\vdash\bf{SD}\ \ {\rm and} \ \ \bf{TP}\vdash\bf{SP}.$$

A brief explanation of our idea, in general, is that, for a theory $\mathcal{T}$, we define its subtheory, i.e. its fragment, in order to understand it better. Namely, if a part of $\mathcal{T}$ is based on axioms $A$ and $B$, and if we are familiar with some simpler statements $A^\prime$ and $B^\prime$, then, in case when $A\vdash A^\prime$ and $B\vdash B^\prime$, we can consider a subtheory $\mathcal{T^\prime}$ based on axioms $A^\prime$ and $B^\prime$, instead of $A$ and $B$. Obviously, $\mathcal{T^\prime}\subseteq\mathcal{T}$ and if $\mathcal{T^\prime}$ is inconsistent, then $\mathcal{T}$ will be inconsistent, but not conversely.

\vskip3truemm
\S 5. \bf Arrow--Sen Theory Simplified. \rm The formal language of symbolic logic presents an adequate framework for refinement of causal connections between intuitive and unclear concepts and conditions formulated initially in a natural language. For instance, a dictatorship condition, as formally introduced by Arrow, helped the author of this article many times to recognise a dictator in real life.

Let $V$ and $X$ be finite sets of individuals and alternatives, respectively, and $P_i$ and $P$ are individual and social strict preference relations on $X$. We define Arrow's dictatorship condition ${\bf{SD}}$, the Pareto property ${\bf{SP}}$, Sen's liberalism axiom ${\bf{SL}}$ and vetoer condition ${\bf{SV}}$, as considered by Fishburn (see [FISH]), all in style of Arrow--Sen social choice theory, but simplified:

${\bf{SD}}$: $(\exists i\in V)(\forall x,y\in X)(xP_iy\to xPy)$

${\bf{SP}}$: $(\forall x,y\in X)((\forall i\in V)xP_iy\to xPy)$

${\bf{SL}}$: $(\forall i\in V)(\exists x,y\in X)(\neg xI_iy\wedge(xP_iy\to xPy)\wedge(yP_ix\to yPx))$

${\bf{SV}}$: $(\exists i\in V)(\forall x,y\in X)(\neg xI_iy\wedge(xP_iy\to\neg yPx))$

Let us emphasize that in ${\bf{SL}}$ and ${\bf{SV}}$ we suppose that $\neg xI_iy$ instead of $x\neq y$, as given in original formulations. It means that, for $\neg xI_iy$ we have $xP_iy\vee yP_ix$, for all alternatives $x$ and $y$ under consideration, and that $\neg xI_iy\to x\neq y$, but not conversely.

Mas--Colell and Sonnenschein [MS] introduced a weak dictatorship condition whose simplified version ${\bf{SWD}}$ looks as follows: $(\exists i\in V)(\forall x,y\in X)(xP_iy\to xRy)$, slightly inconvenient condition because it mixes strict and weak preferences. Here we will prove that conditions ${\bf{SV}}$ and ${\bf{SWD}}$ are logically equivalent. Namely, if we suppose ${\bf{SV}}$, then from $\neg yPx$, bearing in mind that $\neg yPx\leftrightarrow\neg yRx\vee xRy$, and that $R$ is linear, we infer $xRy$, i.e. that ${\bf{SV}}\vdash{\bf{SWD}}$. Conversely, if ${\bf{SWD}}$, then from $xRy$ and $xRy\leftrightarrow xPy\vee xIy$, bearing in mind that $P$ is asymmetric and that $\neg xIy$, we conclude ${\bf{SWD}}\vdash{\bf{SV}}$. This reasoning can be a good exercise for students.

Now we follow Sen's idea of the proof of Arrow's theorem based on the {\it Field--Expansion Lemma} and the {\it Group--Contraction Lemma}. First let us bring to mind the notion of a {\it decisive group} $G$, over two particular alternatives $a,b\in X$, as any nonempty subset of the set $V$ of all individuals, $\emptyset\neq G\subseteq V$, with the property that $((\forall i\in G)aP_ib\to aPb)\wedge((\forall i\in G)bP_ia\to bPa)$. A group $G$ is said to be {\it decisive} if it is decisive over any two alternatives $x,y\in X$, i.e. $(\forall x,y\in X)((\forall i\in G)xP_iy\to xPy)$. A subgroup $F$ of a group $G$ will be any non--empty subset of $G$, $\emptyset\neq F\subseteq G$.

Let us note that the statement that the set $V$ of all individuals is decisive, $V=G$, for a decisive group $G$, is logically equivalent to the Pareto rule ${\bf{SP}}$. On the other hand, the fact that a group consisting of just one element is decisive, will be equivalent to dictatorship axiom ${\bf{SD}}$. So, through the notion of a decisive group this two principles, the Pareto rule ${\bf{SP}}$ and dictatorship ${\bf{SD}}$, are connected immediately.

Two lemmata can be formulated and proved, in accordance with Sen's combinatorial approach (see [Sb], [Sc], [Se] or [BBd]), as follows:

\vskip3truemm
\bf Field--Expansion Lemma. \it If a group is decisive over any two particular alternatives, then it is decisive.\rm\vskip3truemm

Proof. Let $G$ be a decisive group over two particular alternatives $a$ and $b$, and let $x$ and $y$ be any two alternatives, where all four alternatives are mutually distinct. Let us suppose that $\forall i(i\in G\to xP_ia\wedge aP_ib\wedge bP_iy)$, where, by transitivity, we also have $\forall i(i\in G\to xP_iy)$, and $\forall i(i\in V\setminus G\to xP_ia\wedge bP_iy)$, which is possible by the unrestricted domain meta--axiom. It means that: $\forall i\in G(aP_ib)$, $\forall i\in V(xP_ia)$ and $\forall i\in V(bP_iy)$, wherefrom, by decisiveness of $G$ and by the Pareto rule twice, respectively, we draw $aPb$, $xPa$ and $bPy$, and, then, by transitivity, finally, we conclude $xPy$, i.e. that $G$ is decisive over any two alternatives of $X$. Let us emphasize that the meta--axiom of independence of irrelevant alternatives was used here as well. The proof in which the alternatives are not all mutually distinct is quite similar.\qed

\vskip3truemm
\bf Group--Contraction Lemma. \it If a group consisting of at least two individuals is decisive, then there is its proper subgroup which is decisive.\rm\vskip3truemm

Proof. Let us suppose that a decisive group $G$ consists of at least two alternatives, and that its subgroups $E$ and $F$ presents its partition, meaning that $E$ and $F$ are nonempty, $E\cup F=G$ and $E\cap F=\emptyset$. Suppose also that, for any alternatives $x,y$ and $z$, we have $\forall i(i\in E\to xP_iy\wedge xP_iz)$ and $\forall i(i\in F\to xP_iy\wedge zP_iy)$. Similarly as in the bisection interval method, we have the following two possibilities: $xPz$, when $E$ would be decisive over $x$ and $z$, and, if $E$ is not decisive over $x$ and $z$, then $zRx$, where $R$ is a weak version of a strict preference relation $P$. Bearing in mind that $G$ is decisive over $x$ and $y$, and $zRx\wedge xPy\to zPy$, we infer $zPy$. But, only for members $i\in F$, we have $zP_iy$, meaning that, by our hypothesis, $F$ is decisive over $z$ and $y$. Consequently, either $E$ or $F$ must be a decisive group, that means: \it if a group of more than one person is decisive, then so is some its proper subgroup.\qed

\vskip3truemm
\bf Theorem. \rm (Counterpart of Arrow's impossibility of a Paretian non--dictator) \it $\bf{SP}\vdash\bf{SD}$\rm\vskip3truemm

Proof. Immediately, by the Pareto rule, the group of all individuals is decisive. Since it is finite, by successive partitioning, and applying the two above lemmata, each time picking its decisive part, we arrive at a decisive individual, a dictator.\qed

By following Sen's original proof (see [Sb]), or its slight modification as given in [BBb], we also can prove a counterpart of Sen's impossibility theorem:

\vskip3truemm
\bf Theorem. \rm (Counterpart of Sen's impossibility of a Paretian liberal) \it $\bf{SP},\bf{SL}\vdash$\rm\vskip3truemm

Proof. First we will show that simplified liberalism axiom $\bf{SL}$ is equivalent to the following condition
$${\bf{SL^\prime}}: \ \ (\forall i\in V)(\exists x,y)(xP_iy\wedge xPy)$$
Let us consider quantifier free parts of $\bf{SL}$ and $\bf{SL^\prime}$:
\begin{equation}\label{eq1} (xP_iy\vee yP_ix)\wedge(xP_iy\to xPy)\wedge(yP_ix\to yPx)\end{equation}
and
\begin{equation}\label{eq2}  \ \ \ (xP_iy\wedge xPy)\vee(yP_ix\wedge yPx)\end{equation}
respectively. If from the first parentheses of $(1)$ we have $xP_iy$, then, by {\it modus ponens} and the second parentheses of $(1)$, we can derive the first parentheses of $(2)$; but, from $yP_ix$, in a similar way, we infer the second parentheses of $(2)$. Consequently, ${\bf{SL}}\vdash{\bf{SL^\prime}}$. Conversely, the first parentheses of $(2)$, $xP_iy$ and $xPy$, obviously, enable us to derive the first and second parentheses of $(1)$, but, additionally, from $xP_iy$, by asymmetry, we have $\neg yP_ix$, i.e. $\neg yP_ix\vee yPx$, which are the third parentheses of $(1)$. Consequently, ${\bf{SL^\prime}}\vdash{\bf{SL}}$.

Now, let us prove our theorem. Suppose that the set $V$ of all individuals consists of $n(\geq2)$ persons and that, for the alternatives $x_1,\dots,x_{n}\in X$ and $y_1,\dots,y_{n}\in Y$, the particular cases of the liberalism axiom
$$(x_{i}P_{i}y_i\wedge x_{i}Py_i)\vee(y_{i}P_{i}x_i\wedge y_{i}Px_i)$$
hold, for all $i$ $(1\leq i\leq n)$. Let us analyze the preferences of the first two individuals only, having in mind that, for alternatives $x_1,y_1,x_2$ and $y_2$, $x_1Py_1\wedge x_2Py_2$, $x_1Py_1\wedge y_2Px_2$, $y_1Px_1\wedge x_2Py_2$ or $y_1Px_1\wedge y_2Px_2$holds. For the first combination $x_1Py_1\wedge x_2Py_2$, in the case when $y_1=x_2$ or $x_1=y_2$, it is possible to suppose that $y_2P_ix_1$, for all $i$ $(1\leq i\leq n)$, or  $y_1P_ix_2$, for all $i$ $(1\leq i\leq n)$, respectively, wherefrom, by the Pareto rule, we can infer $y_2Px_1$ or $y_1Px_2$, violating that $\neg x_2Px_2$ or $\neg x_1Px_1$, by transitivity, in each case. In a similar way other combinations can be discussed. Supposing again $x_1Py_1\wedge x_2Py_2$, if $y_1\neq x_2\wedge x_1\neq y_2$, then it is possible to suppose that, for all $i$ $(1\leq i\leq n)$,  $y_2P_ix_1\wedge y_1P_ix_2$, wherefrom, by the Pareto rule, we can infer $y_2Px_1\wedge y_1Px_2$, violating that $\neg y_2Py_2$, by transitivity, in each case, as well. Note that, in both cases we used the unrestricted domain meta--axiom implicitly.\qed

In the above proof we followed, essentially, the spirit of Sen's original proof (see [Sb] or [Se]), based on the minimal liberalism argument, meaning that there are at least two persons decisive over two existing pairs of distinct alternatives.

In order to introduce vetoing condition easily into our consideration, first we prove:\vskip3truemm

\bf Lemma. \rm (Counterpart of Fishburn's note on a dictatorial non--vetoer) \it $\bf{SD}\vdash\bf{SV}$\rm\vskip3truemm

Proof. \rm Obviously, from ${\bf{SD}}$: $xP_iy\to xPy$, bearing in mind that $P$ is asymmetric, $(\forall x,y\in X)(xPy\to\neg yPx)$, we conclude ${\bf{SV}}$: $xP_iy\to\neg yPx$, i.e. that ${\bf{SD}}\vdash{\bf{SV}}$, but not conversely.\qed

As vetoer and weak dictatorship conditions are equivalent, we can formulate the following statement:

\vskip3truemm
\bf Theorem. \rm (Counterpart of Mas--Colell---Sonnenschein's impossibility of a Paretian non--vetoer, or impossibility of a Paretian non--weak dictator) \it $\bf{SP}\vdash\bf{SV}$\rm\vskip3truemm

This new context enables us to prove and present a set of impossibility results more easily. Moreover, we can establish some new statements.

\vskip3truemm
\bf Lemma. \rm (Counterpart of Chichilnisky's impossibility of a non--Paretian dictator) \it $\bf{SD}\vdash\bf{SP}$\rm\vskip3truemm

Proof. \rm By following the proof presented in [BBe], from $xP_iy\to xPy$, with the help of the weakening antecedent rule $n$--times, by which from $p\to q$ one can always derive $p\wedge r\to q$, we infer $\bigwedge_{1\leq j\leq n}xP_jy\to xPy$, i.e. $(\forall j\in V)xP_jy\to xPy$. It means that $\bf{SD}\vdash\bf{SP}$.\qed

An immediate consequence is that: $\bf{SD}\vdash\bf{SP}$. Bearing in mind that Arrow's impossibility of a Paretian non--dictator can also be expressed as $\bf{SP}\vdash\bf{SD}$, we obtain a logical equivalence between the Pareto property and dictatorship. This result was obtained by Chichilnisky [CH] in a topological context.

\vskip3truemm
\bf Lemma. \rm (Counterpart of Impossibility of a liberal dictatorship) \it $\bf{SD},\bf{SL}\vdash$\rm\vskip3truemm

Proof. \rm Suppose that $i$ is a dictator and that his personal preference regarding particular alternatives $x$ and $y$ is $xP_iy$, having as its consequence the social preference $xPy$. On the other side, let $j$ be a person, $j\neq i$, possessing liberty power over the same two alternatives $x$ and $y$ such that $yP_jx$ having as a consequence the social preference $yPx$, wherefrom, by asymmetry of $P$, we infer $\neg xPy$, which is inconsistent with the previous dictatorial conclusion $xPy$.\qed

The example that follows is a trivial and expectable statement related to liberalism and veto power.

\vskip3truemm
\bf Lemma. \rm (Counterpart of Impossibility of a liberal vetoing) \it $\bf{SL},\bf{SV}\vdash$\rm\vskip3truemm

Proof. \rm Let $i$ be a person with preference $yP_iy$ possessing liberty power over two alternatives $x$ and $y$ such that its consequence is social preference $yPx$. Suppose that a vetoer $j$, $j\neq i$, defines his personal preference $xP_iy$ with consequence $\neg yPx$, that makes a system inconsistent due to the fact that $yPx$.\qed

We are able to give a simple argumentation for Mas--Colell---Sonnenschein's result:

\vskip3truemm
\bf Lemma. \rm (Counterpart of Impossibility of a Paretian non--vetoer) \it $\bf{SP}\vdash\bf{SV}$\rm\vskip3truemm

Proof. \rm This conclusion, $\bf{SP}\vdash\bf{SV}$, can be inferred directly from Arrow's 'impossibility of a Paretian non--dictator', $\bf{SP}\vdash\bf{SD}$, and Fishburn's 'impossibility of a dictatorial non--vetoer', $\bf{SD}\vdash\bf{SV}$, i.e. that each dictator is a vetoer.\qed%

An immediate consequence is that each Paretian society has a vetoer: $\bf{SP}\vdash\bf{SV}$. The similar conclusion can be inferred from Arrow's impossibility of a Paretian non--dictator, $\bf{SP}\vdash\bf{SD}$, and Fishburn's impossibility of a dictatorial non--vetoer, $\bf{SD}\vdash\bf{SV}$.

From the previous two lemmata, Impossibility of a Paretian non--vetoer, $\bf{SP}\vdash\bf{SV}$, and impossibility of a liberal vetoing, $\bf{SL},\bf{SV}\vdash$, we derive immediately, by using the cut rule, Sen's theorem:

\vskip3truemm
\bf Corollary. \rm (Counterpart of Impossibility of a Paretian liberal) \it $\bf{SP},\bf{SL}\vdash$

\vskip3truemm\rm Let us consider now a simplified version of the strong Pareto rule

${\bf{SSP}}$: $(\forall x,y\in X)((\forall i\in V)xR_iy\wedge(\exists i\in V)xP_iy\to xPy)$

\noindent and a simplification of the strong dictatorship

${\bf{SSD}}$: $(\exists i\in V)((\exists x,y\in X)(xP_iy\to xPy)\wedge(\exists x,y\in X)(xR_iy\to xRy))$

\noindent as defined in [Sa]. Obviously, ${\bf{SSP}\vdash\bf{SP}}$ and ${\bf{SD}\vdash\bf{SSD}}$, but also ${\bf{SP}\not\vdash\bf{SSP}}$ and ${\bf{SSD}\not\vdash\bf{SD}}$.

Although it sounds very simple, herefrom, by immediate turnstile formalism and the above conclusions, we can infer the following additional impossibilities:

\vskip3truemm
\bf Corollary. \it (a) $\bf{SSP},\bf{SD}\vdash$;

\ (b) $\bf{SSP},\bf{SL}\vdash$;

\ (c) $\bf{SP}\vdash\bf{SSD}$;

\ (d) $\bf{SSP}\vdash\bf{SSD}$.

\vskip3truemm\rm This could be a good exercise for students.

\vskip3truemm
\S 6. \bf Concluding Remarks. \rm Simplifications presented here may be divided into at least three lines: changing the status of two higher--order conditions, 'unrestricted domain' and 'the independence of irrelevant alternatives', abstracting the notions of a social welfare function and a social decision function, and, finally, substituting some easier axioms, such as dictatorship, vetoing, liberalism and the Pareto property, by their first--order consequences. The basic result is that after these simplifications, we obtain a fragment of traditional Arrow--Sen theory in which we can also prove well--known interdeducibilities  and impossibilities, including Arrow's and Sen's theorems. Further, due to the finiteness of sets of individuals $V$ and alternatives and $X$, quantifiers can be replaced by finite conjunctions and disjunctions. This can be a good starting point in crossing from the first--order language to the propositional language, extended by a finite list of binary relations defined over the set of propositional formulae. In this spirit, the multiple von Wright preference logic (see [BBd] and [BBe]) makes a suitable framework for mechanical and formal dealing with preferences and alternatives. This extension of classical propositional calculus, with our conjecture that this is a decidable logic, presents a good theoretical basis for programming decision making support systems as part of the wider artificial intelligence project.

Finally, let us consider an example concerning famous Arrow's impossibility theorem. This theorem can be formulated as $\bf{TP}\vdash\bf{TD}$, assuming that conditions of 'unrestricted domain' and 'the independence of irrelevant alternatives' hold, in original Arrow's theory, while its analogue, a similar statement, in a new simplified context, is the following one: $\bf{SP}\vdash\bf{SD}$. Let us emphasize that neither counterpart $\bf{SP}\vdash\bf{SD}$ implies original Arrow's theorem $\bf{TP}\vdash\bf{TD}$, nor vice versa. Consequently, these two statements can be considered as two roughly connected facts in two parallel worlds. Similarly, we can present a counterpart of Chichilnisky's original theorem [C82], 'impossibility of a non--Paretian dictator': $\bf{TD}\vdash\bf{TP}$, and its counterpart in our simplified context: $\bf{SD}\vdash\bf{SP}$, asserting again that there is no immediate formal logical connection between these two statements. But, on the other side, bearing in mind that $\bf{TP}\vdash\bf{SP}$ and $\bf{TL}\vdash\bf{SL}$, we can directly derive well--known Sen's 'impossibility of a Paretian liberal': $\bf{TP},\bf{TL}\vdash$, from its simplified version $\bf{SP},\bf{SL}\vdash$, meaning that the axioms $\bf{SP}$ and $\bf{SL}$, and, consequently, the axioms $\bf{TP}$ and $\bf{TL}$, when appear together, make a theory inconsistent.

We consider that the value of this simplified approach, focusing on 'simpler' (first--order) conditions, is in giving an opportunity to a wider circle of readers to understand the basic ideas, results and spirit of traditional Social Choice Theory better.

\vskip2truemm 

\bf References: \rm

\vskip2truemm
\noindent [AHW] T. Agotnes, W. van der Hoek, M. J. Wooldridge, \it On the logic of preference and judgment aggregation, \bf Auton. Agents Multi Agent Syst. \rm 22(1) (2011), pp. 4--30.

\vskip2truemm
\rm\noindent [AR] K. J. Arrow, \bf Social Choice and Individual Values, \rm John Wiley, New York, 1963.

\vskip2truemm
\rm\noindent [AS] K. J. Arrow, A. Sen, \bf Handbook of Social Choice and Welfare, \rm vol 1, Gulf Professional Publishing, Houston, 2002.

\vskip2truemm
\rm\noindent [ASS] K. J. Arrow, A. Sen, K. Suzumura, \bf Handbook of Social Choice and Welfare, \rm vol 2, Elsevier, Amsterdam, 2010.


\vskip2truemm
\noindent [BL] M. Balinski, R. Laraki, \bf Majority Judgement --- Measuring, Ranking and Electing, \rm MIT Press, Cambridge, MA, 2010.

\vskip2truemm
\noindent [BBa]
B. Bori\v ci\'c, {\it Logical and historical determination of impossibility theorems by Arrow and Sen}, {\bf Economic Annals} LI 172 (2007), pp. 7--20.

\vskip2truemm
\noindent [BBb]
B. Bori\v ci\'c, {\it Dictatorship, liberalism and the Pareto rule: possible and impossible,} {\bf Economic Annals} LIV 181 (2009), pp. 45--54.

\vskip2truemm
\noindent [BBc]
B. Bori\v ci\'c, {\bf Logic and Proof}, Faculty of Economics, University of Belgrade, Belgrade, 2011. (Zbl 1214.03043; MR 2011k:03001)

\vskip2truemm
\noindent [BBd]
B. Bori\v ci\'c, {\it Multiple von Wright’s preference logic and social choice theory,} {\bf The Bulletin of Symbolic Logic} 20, 2014, pp. 224--225.

\vskip2truemm
\noindent [BBe]
B. Bori\v ci\'c, {\it Impossibility theorems in multiple von Wright's preference logic,} {\bf Economic Annals}, LIX 201 (2014), pp. 69--84.

\vskip2truemm
\noindent [BBf]
B. Bori\v ci\'c, {\it A note on dictatorship, liberalism and the Pareto rule,} {\bf Economic Annals}, LXVIII 238 (2023), pp. 115--119.


\vskip2truemm
\noindent [BS]
B. Bori\v ci\'c, M. Sre\'ckovi\'c, {\it Vetoing: social, logical and mathematical aspects}, {\bf Mathematics for Social Sciences and Arts, Algebraic Modeling} (Eds. M. N. Hounkonnou et al.), Springer, 2024, pp. 101--124.


\vskip2truemm
\noindent [BUCa]
J. M. Buchanan, {\it Social choice, democracy, and free markets}, {\bf Journal of Political Economy} 62(2)(1954), pp. 114--123.

\vskip2truemm
\noindent [BUCb]
J. M. Buchanan, {\it Individual choice in voting and the market}, {\bf Journal of Political Economy} 62(3)(1954), pp. 334--343.

\vskip2truemm
\noindent [CH] G. Chichilnisky, \it The topological equivalence of the Pareto condition and the existence of a dictator,
\bf Journal of Mathematical Economics \rm 9 (1982) pp. 223--233.


\vskip2truemm\noindent [DS]
H. de Swart, \bf Philosophical and Mathematical Logic, \rm Springer Undergraduate Texts in Philosophy, Springer, Berlin, 2018.

\vskip2truemm
\noindent [FISH]
P. C. Fishburn, {\bf The Theory of Social Choice,} Princeton University Press, Princeton, 1973.

\vskip2truemm
\noindent [GB]
D. M. Gabbay, \bf Fibring Logics, \rm Clarendon Press, Oxford, 1999.\rm




\vskip2truemm
\noindent [KEL]
J. S. Kelly, {\bf Arrow Impossibility Theorems,} Academic Press, London, 1978.

\vskip2truemm
\noindent [MS]
A. Mas--Colell, H. Sonnenschein, {\it General possibility theorems for group decisions,} {\bf The Review of Economic Studies} 39 (1972), pp. 185--192.

\vskip2truemm
\noindent [MSK]
E. Maskin, {\it A modified version of Arrow's $\bf{IIA}$ condition,} {\bf Social Choice and Welfare} 54 (2020), pp. 203--209.

\vskip2truemm
\noindent [NIP] T. Nipkow, \it Social choice theory in HOL, \bf J. Autom. Reasoning \rm 43(3) (2009), pp. 289--304.

\vskip2truemm
\noindent [PY] E. Pacuit, F. Yang, \it Dependence and independence in social choice: Arrow's theorem, \rm in S. Abramsky, J. Kontinen, H. Vollmer, J. Väänänen, eds, \bf Dependence Logic: Theory and Applications, Progress in Computer Science and Applied Logic, \rm Birkhauser, June 2016, pp. 235--260.

\vskip2truemm
\noindent [MP] M. Pauly, \it On the role of language in social choice theory, \bf Synthese \rm 163(2) (2008), pp. 227--243.

\vskip2truemm
\noindent [PER] T. Perkov, \it Natural deduction for modal logic of judgment aggregation, \bf Journal of Logic, Language and Information \rm 25(3--4) (2016),pp. 335--354.

\vskip2truemm
\noindent [POR] D. Porello, \it Logics for modelling collective attitudes, \bf Fundam. Inform. \rm 158(1-3) (2018), pp. 239--275.

\vskip2truemm
\noindent [RTH] R. Routhley, \it Repairing proofs of Arrow's general impossibility theorem and enlarging the
scope of the theorem, \bf Notre Dame Journal of Formal Logic \rm 20 (1979), pp. 879--890.

\vskip2truemm
\noindent [SCH] N. J. Schofield, \bf Social Choice and Democracy, \rm Springer--Verlag, Berlin, 1985.

\vskip2truemm
\noindent [Sa] A. K. Sen, \it Quasi--Transitivity, Rational Choice and Collective Decisions, \bf The Review of Economic Studies \rm 36 (3) (1969), pp.
381--393.

\vskip2truemm
\noindent [Sb] A. K. Sen, \it The impossibility of a Paretian liberal, \bf Journal of Political Economy \rm 78 (1) (1970), pp.
152--157.

\vskip2truemm
\noindent [Sc] A. K. Sen, \bf Collective Choice and Social Welfare, \rm Holden--Day, San Francisco, (1970) (Fourth
edition: Elsevier, Amsterdam, 1995).

\vskip2truemm
\noindent [Sd] A. K. Sen, \it Internal Consistency of Choice, \bf Econometrica \rm 61 No. 3. (1993), pp. 495--521.

\vskip2truemm
\noindent [Se] A. K. Sen, \it Rationality and social choice, \bf American Economic Review \rm 85 (1995), pp. 1--24.


\vskip2truemm
\noindent [SRE]
M. Sre\'ckovi\'c, {\it Measuring the veto power,} {\bf Ekonomske ideje i praksa} 25, 2017, pp. 89--96.

\vskip2truemm
\noindent [SZ]
K. Suzumura, {\it Reflections on Arrow's research program of social choice theory,} {\bf Social Choice and Welfare} 54 (2020), pp. 219--235.

\vskip2truemm
\noindent [TAK]
G. Takeuti, {\bf Proof Theory}, North--Holland Publishing Company, Amsterdam, 1975  (Second edition: Dover Publications, 2013).

\vskip2truemm
\noindent [TL] P. Tang, F. Lin, \it Computer--aided proofs of Arrow's and other impossibility theorems, \bf Artif. Intell. \rm  173(11) (2009), pp. 1041--1053.

\vskip2truemm
\noindent [TAY]
A. D. Taylor, {\bf Mathematics and Politics --- Strategy, Voting, Power and Proof,} Springer--Verlag, Berlin, 1995.

\vskip2truemm
\noindent [THW] N. Troquard, W. van der Hoek, M. J. Wooldridge, \it Reasoning about social choice functions, \bf J. Philosophical Logic \rm 40(4) (2011) pp. 473--498.

\vskip2truemm
\noindent [DvD]
D. van Dalen, {\bf Logic and Structure}, Spriner, Berlin, 1980 (Fifth edition 2013).

\vskip2truemm
\noindent [VWa]
G. H. von Wright, \bf The Logic of Preference, \rm Edinburgh University Press,
Edinburgh, 1963.

\vskip2truemm
\noindent [VWb]
G. H. von Wright, \it The logic of preference reconsidered, \bf Theory and Decision \rm 3 (1972) pp. 140--169, or \bf Philosophical Logic --- Philosophical Papers, Volume II, \rm Cornell University Press, New York, 1983, pp. 67--91.

\vskip2truemm
\noindent [WIN]
E. Winter, {\it Voting and Vetoing,} {\bf The American Political Science Review} 90 (1996), pp. 813--823.\vskip2truemm
\end{document}